\documentclass[12pt]{article}
\usepackage{amsmath}
\usepackage{amssymb}
\usepackage{amsfonts}
\usepackage{eucal}
\usepackage[usenames]{color}
\usepackage{graphicx}
\numberwithin{equation}{section}
\newtheorem{Theorem}{Theorem}[section]
\newtheorem{Proposition}[Theorem]{Proposition}

\newtheorem{Remark}[Theorem]{Remark}
\begin{document}

\title{Inverse first-passage problems of a  diffusion with  resetting}

\author{Mario Abundo\thanks{Dipartimento di Matematica, Universit\`a  ``Tor Vergata'', via della Ricerca Scientifica, I-00133 Rome, Italy.
E-mail: \tt{abundo@mat.uniroma2.it}}
}
\date{}
\maketitle

\begin{abstract}
We address some inverse problems for the first-passage place and the first-passage time of a one-dimensional diffusion process $\mathcal X(t)$  with stochastic resetting,
starting from
an initial position $\mathcal X(0)= \eta ;$ this type of diffusion $\mathcal X(t)$ is characterized by the fact that a reset to the position $x_R $ can occur
according to a homogeneous Poisson process with rate $r>0.$
As regards the  inverse first-passage place problem,
for random  $\eta \in (0,b), \ b < + \infty$ (and fixed $r$ and $x_R \in (0,b))$, let  $\tau_{0,b}$ be the first time at which $\mathcal X(t)$ exits the interval $(0,b),$
and $\pi _0 = P(\mathcal X(\tau_{0,b}) =  0)$ the probability of exit from the left end of $(0,b);$
given a probability $q \in (0,1),$
the  inverse first-passage place problem consists in finding the density $g$ of  $\eta ,$ if it
exists, such that $\pi _0 = q.$ Concerning the inverse first-passage time problem, for random $\eta \in (0, + \infty)$ (and fixed $r$ and $x_R >0)$, let $\tau$ be the first-passage time of $\mathcal X(t)$ through zero;
for a given distribution function $F(t)$ on the positive real axis, the inverse first-passage time  problem  consists in finding the density $g$ of  $\eta,$  if it exists, such that $P(\tau \le t ) = F(t), \ t >0.$
In addition to the case of random initial position $\eta,$ we also study the case  when the initial position  $\eta$ and the resetting rate $r$ are fixed, whereas the reset position $x_R$ is random.
For all types of inverse problems considered, several explicit examples of solutions are reported.
\end{abstract}

\par{\bf Subjclass [2020]:} {Primary 60J60, 60H05; Secondary 60H10}\par
{\bf Keywords:} {Diffusion with resetting, first-passage time, first-passage place.}



\section{Introduction}
This paper has to do in part with the article \cite{abundo:FPA2023} that concerns the first-passage area of a one-dimensional diffusion process with stochastic resetting $\mathcal X(t);$ this type of process can describe the temporal evolution of a dynamical
system which evolves starting from a given initial value, but only for a random period of
time, after that the dynamics gets renewed and starts afresh from a reset value, continuing
its evolution until the next resetting event occurs, and so on.
The process $\mathcal X(t)$ is described
precisely below. \par
Let $X(t)$ be a one-dimensional temporally homogeneous diffusion process,
driven by the SDE:
\begin{equation} \label{diffusion}
dX(t)=  \mu(X(t)) dt + \sigma (X(t)) d{W_t} ,
\end{equation}
and starting from
an initial position $X(0)=\eta $ (fixed or random),
where ${W_t}$ is a standard Brownian motion (BM) and the drift $\mu (\cdot)$ and diffusion coefficient $\sigma (\cdot)$ are regular enough functions, such that there exists a unique strong solution of the
SDE \eqref{diffusion} (see e.g. \cite{klebaner}).
\par\noindent
From $X(t)$ we construct a new process $\mathcal X(t),$ as follows.
We suppose that resetting events can occur according to a homogeneous Poisson process with rate $r>0.$
Until the first resetting event the process $\mathcal X(t)$ coincides with $X(t)$ and it evolves according to \eqref{diffusion} with $\mathcal X(0)=X(0)=\eta ;$ when the reset occurs,
$\mathcal X(t)$ is set instantly to a position $x_R .$ After that, $\mathcal X(t)$ evolves again according to \eqref{diffusion} starting afresh (independently of the past history) from $x_R,$ until the next resetting event occurs, and so on. The inter-resetting times turn out to be independent and exponentially distributed random variables with parameter $r.$
In other words, in any time interval $(t, t+ \Delta t),$ with $\Delta t \rightarrow 0 ^+, $ the process can pass from $\mathcal X(t)$  to the position $x_R$ with probability $r \Delta t  + o( \Delta t),$ or it can continue its evolution according to \eqref{diffusion} with probability $1- r \Delta t + o( \Delta t ).$
\par\noindent
The process $\mathcal X(t)$ so obtained is called diffusion process with stochastic resetting;
it has some analogies with the process considered in \cite{dicre:03}, where it was studied a M/M/1 queue with catastrophes and its  continuous approximation, namely a Wiener process subject to randomly occurring jumps at a given rate $\xi,$ each jump making the process instantly obtain the state $0.$
 Thus, the process considered in \cite{dicre:03} can be viewed as a Wiener process with resetting,
 in which a reset to the position $x_R =0$ is done, according to a homogeneous Poisson process with rate $r= \xi.$ \par
For any $C^2$ function $f(x),$ the infinitesimal generator of $\mathcal X(t)$ is given by (see e.g. \cite{abundo:FPA2023}):
\begin{equation} \label{generator}
\mathcal L f(x) = \frac 1 2 \sigma ^2(x) f''(x) + \mu (x) f'(x) +r (f(x_R) -f(x)) \equiv L f(x) +r (f(x_R) -f(x)) ,
\end{equation}
where $Lf(x)= \frac 1 2 \sigma ^2(x) f''(x) + \mu (x) f'(x)$ is the  ``diffusion part'' of the generator, i.e. that concerning the diffusion process $X(t).$ \par\noindent
\par
In this paper, we address some inverse problems for the first-passage place and the first-passage time of a diffusion with stochastic resetting $\mathcal X(t).$ Two cases are considered: in case I, the initial position
$\eta = \mathcal X(0)$
is supposed to be random and independent of $\mathcal X(t),$ whereas the resetting rate $r$ and the reset position $x_R$ are fixed; in  case II, the initial position $\eta = \mathcal X(0)$ and $r$ are fixed, whereas  $x_R$ is random  and independent of $\mathcal X(t).$ \par
\underline {Case I}.  As regards the  inverse first-passage place (IFPP) problem, we suppose that $\mathcal X(t)$ is a diffusion with resetting, and  we take the initial position $\eta = \mathcal X(0) $ randomly distributed in a bounded interval $(0,b)$ and independent of $\mathcal X(t)$  (the reset rate $r>0$ and the reset position $x_R \in (0,b)$ are fixed).
We suppose that the first-exit time (FET) of $\mathcal X(t)$ from the interval $(0,b):$
\begin{equation}
\tau _ {0, b} = \inf \{t \ge 0: \mathcal X(t) \notin (0, b) \},
\end{equation}
is finite with probability one, and we denote by
$\mathcal X( \tau _{0,b})$ the first-passage place of $\mathcal X(t)$ at time $\tau _{0,b};$ note that, by assumption,
one has $\mathcal X( \tau _{0, b}) =0 ,$ or $\mathcal X(\tau _{0, b}) = b.$ Let  $\pi _0 = P ( \mathcal X(\tau _{0, b}) =0)$ be the probability that the process $\mathcal X(t)$ first exits the interval $(0,b)$
from the left end, and
$\pi _b = 1- \pi _0 = P ( \mathcal X(\tau _{0, b} ) = b )$ the probability  that $\mathcal X(t)$  first exits from the right. \par\noindent
Then, for fixed $b, \ r>0$ and $x_R \in (0,b),$
we consider the following IFPP problem for $\mathcal X(t)$ in $(0,b):$
\begin{equation} \label{IFPPetarandom}
\begin{aligned}
given \ a \ probability \ q \in (0,1), \ find \ the \ density \\
g \ of \ \eta \in (0,b) \ (if \ it \ exists), \ so \ that \ \pi_0 = q .
\end{aligned}
\end{equation}
\medskip

\noindent The {density} $g$ is called a solution to the IFPP
problem \eqref{IFPPetarandom}. In fact, the solution to such a problem, if it exists, is not necessarily unique (see Remark 2.2, and also  \cite{abundo:CC22}, \cite{abundo:saa20IFPP}, {\bf {\cite{lefeb:22}}}
 for diffusions without resetting).
Of course, we can also admit $\eta$ to be a discrete random variable taking values in a subset of $[0,b];$ in this case $g(x)$ turns out to be a discrete probability density.
\par\noindent
For general jump-diffusion processes without resetting, even though  several papers on first-passage times are available (see e.g. \cite{abundo:mcap10}, \cite{abundo:pms00}, \cite{kouwang:03}, \cite{tuckwell:76}, and references therein), few results are known for first-passage places;
the direct first-passage place problem was studied in \cite{kouwang:03} and \cite{lefeb:19}, where
equations for the moments of first-passage places were established;
the inverse first-passage place problem was studied  in \cite{abundo:CC22}, \cite{abundo:saa20IFPP}, { \cite{lefeb:22}}. \par
As concerns the inverse first-passage time (IFPT) problem for a diffusion with resetting $\mathcal X(t),$  starting from a random  position $\eta = \mathcal X(0) \in (0, + \infty),$ which is supposed to be independent of  $\mathcal X(t)$ (the reset rate $r>0$ and
the reset position $x_R >0$ are fixed), let $\tau  = \inf \{t >0: \mathcal X(t) =0 \}$ be the first-passage time (FPT) of $\mathcal X(t)$ through zero; then, we consider the following IFPT problem:
\begin{equation} \label{IFPTetarandom}
\begin{aligned}
given \ a \ distribution  \ function \ F(t) \ for \ t>0, \ find \ the \ density \ of \ \eta >0 \\
(if \ it \ exists), \ such \ that \ the \ FPT  \ of \ \mathcal X(t) \ through \ zero \
has \ distribution  \ F(t).
\end{aligned}
\end{equation}
\medskip

\noindent For diffusions without resetting, the analogous IFPT problem  was studied e.g.
in \cite{abundo:saa19}, \cite{abundo:mathematics},
\cite{abu:smj13}, \cite{abundo:stapro12}, \cite{jackson:stapro09}; the corresponding problem for two boundaries was studied in \cite{abundo14b}, \cite{abundo:stapro13},  while the article
\cite{abundo13b}
dealt with the case of a diffusion with jumps (see also \cite{abundo:LNSIM} for a review). \par
In addition to the aforementioned problems,  we study two further types of inverse problems, namely, the IMFPT problem, when the mean of the FPT of $\mathcal X(t)$ through zero
is assigned and one has to find the density of the starting position $\eta \in (0, + \infty),$ and the IMFET problem, when the  mean of the first-exit time (FET) of
$\mathcal X(t)$ from $(0,b)$ is assigned, and one has to find the density of the starting position $\eta \in (0,b).$ \par
\underline {Case II}. We state and study all inverse problems previously considered, also when the initial position $\eta = \mathcal X(0)$ and the reset rate $r$ are fixed, whereas the reset position $x_R$ is random and independent of $\mathcal X(t):$ instead of the density of $\eta ,$ we search for the density of $x_R,$ as a solution to
the inverse problems. \par
Of course, all the inverse problems can be studied in any interval $(a,b),$ or $(a, + \infty);$ we have taken $a=0,$ only for the sake of simplicity.
\par
Note that, in general, the  IFPP, IMFPT, and IMFET problems (or the corresponding ones in case II, see Section 4) can also have more than one solution; instead, the solution to the IFPT problem (or the corresponding one in case II, see Section 4) is unique, provided that its Laplace transform is analytic.
\bigskip

Direct and inverse problems for the first-passage time and the  first-passage place of diffusion processes are worthy of attention, since they
have interesting applications in several applied fields, for instance
in biological modeling for neuronal activity  (see e.g. \cite{lanska:89}, \cite{norisa:85} and the references contained in \cite{abundo:FPA2023}).
They are also relevant in
Mathematical Finance, in particular in credit risk modeling, where
the first-passage time represents a default event of an obligor (see e.g.
\cite{jackson:stapro09}); other applications can be found e.g.
in queuing theory, where the first hitting time to zero can be identified with the
busy period, that is the time until the queue is first empty, and in many other fields (see e.g. the discussion in \cite{abundo:stapro12}).
Finally, for a review concerning functionals of Brownian motion with resetting in physics and computer science, see e.g \cite{maj07}.
\par
At our knowledge,  extensions of inverse problems to diffusions with resetting have not been treated  in the literature, yet; thus, the aim of the present article is just to study these types of problems for one-dimensional
diffusions with resetting $\mathcal X(t)$. We do not deal here with the inverse problem  for $\mathcal X(t),$ analogous to \eqref{IFPTetarandom}, which consists in finding the density of the initial position $\eta \in (0,b)$ (or of the reset position $x_R \in (0,b)),$ in such a way that the first-exit time of  $\mathcal X(t)$ from $(0,b)$ has an assigned distribution; this and the corresponding direct problem will be the subject of a future work.
\par
The paper is
organized as follows: Section 2 and 3 concern case I, namely, when the initial position $\eta = \mathcal X(0)$ is random, whereas the reset rate $r$ and the reset position $x_R$ are fixed; in particular, Section 2 contains the results on the IFPP  problem of  $\mathcal X(t)$;
Section 3 deals with the IFPT problem (subsection 3.1), and  with the IMFPT problem   (subsection 3.2),
and the IMFET problem (subsection 3.3). Section 4 deals with all the inverse problems described above, in case II, that is,  when the initial position
$x = \mathcal X(0) $ and the reset rate $r$ are fixed, whereas the  reset position $x_R$ is random.
For all types of inverse problems considered,
several explicit examples of solutions are reported, mostly concerning Wiener process with resetting (i.e. drifted, or undrifted BM with resetting).
Finally, Section 5 contains conclusions and final remarks.

\section{Case I: The inverse first-passage place (IFPP) problem for random initial position $\eta$ and fixed $r$ and $x_R$}
In this section, we study the IFPP problem \eqref{IFPPetarandom} for a diffusion process with resetting $\mathcal X(t),$
obtained from the underlying diffusion $X(t)$ driven by the SDE \eqref{diffusion},
where the initial position $\eta = \mathcal X(0) $ is supposed to be randomly distributed in a bounded interval $(0,b)$ and independent of  $\mathcal X(t),$ whereas the reset rate $r >0$ and the reset position $x_R \in (0,b)$ are fixed.
If $\tau _{0,b}(x) = \inf \{ t \ge 0: \mathcal X (t) \notin (0,b) | \mathcal X(0)=x \}$  denotes the first-exit time (FET) of $\mathcal X(t)$ from the interval $(0,b),$ under the condition that
$\eta =x \in (0,b),$ we suppose that $\tau _{0,b}(x)$ is finite with probability one;
we recall (see e.g. \cite{abundo:FPA2023}, \cite{abundo:CC22}, \cite{abundo:saa20IFPP}, \cite{abundo:pms00}) that,
 the function $ \pi _0(x) = P(\mathcal X(\tau _{0,b}(x))=0)$
satisfies the problem with boundary conditions:
\begin{equation}
\begin{cases}
\mathcal L f(x)= L f(x) +r (f(x_R) -f(x)) = 0, \ x \in (0,b) \\
f(0)= 1, \  f(b)= 0
\end{cases}
\end{equation}
$( \mathcal L $ and $L$ are defined by \eqref{generator}),
namely:
\begin{equation} \label{eqforpi0}
\begin{cases}
\frac 1 2 \sigma ^2 (x) f'' (x) + \mu (x) f' (x) + r (f(x_R) - f(x))  = 0, \ x \in (0,b) \\
f(0)= 1, \  f(b)= 0,
\end{cases}
\end{equation}
while the mean value of the FET, $E [\tau _{0,b}(x)],$ satisfies:
\begin{equation} \label{eqfortau}
\begin{cases}
L f(x) +r (f(x_R) -f(x)) = -1, \ x \in (0,b) \\
f(0)= f(b)= 0.
\end{cases}
\end{equation}
\bigskip

\noindent
By assuming that the random initial position $\eta \in (0,b)$
possesses a density $g(x),$ we obtain the following result, which is analogous to that holding for diffusions without resetting (see \cite{abundo:saa20IFPP}):
\begin{Proposition} \label{mainpro}
Let $\mathcal X(t)$ be the diffusion process with resetting, described in Section 1, and let be $q \in (0,1);$
with the previous notations,
if a solution $g$ exists to the IFPP problem \eqref{IFPPetarandom}, then the {density} $g$ must satisfy the following equation:
\begin{equation} \label{maineq}
q= \int _0 ^b g(x) \pi_0 (x) dx ,
\end{equation}
where $\pi_0 (x)$ is the solution of the problem \eqref{eqforpi0}.
\end{Proposition}
\hfill  $\Box$
\par\noindent
 {If $g(x) \pi_0(x)$ is continuous in $[0,b],$ then from \eqref{maineq} the mean value theorem implies that there exists $c \in [0,b]$ such that
$g(c)= \frac q {b \pi_0(c)}.$
}
The equation \eqref{maineq} can be written as:
\begin{equation} \label{maineq2}
q= E [\pi_0 (\eta)].
\end{equation}
\begin{Remark}
{\it
For an assigned $q \in (0,1),$ Eq. \eqref{maineq} is an integral equation in the
unknown $g(x).$ Unfortunately, no method is available to solve analytically this equation, so
any possible solution $g$ to the IFPP problem \eqref{IFPPetarandom} must be found by making attempts (see also Remark 2.5 in \cite{abundo:saa20IFPP}).
The IFPP problem \eqref{IFPPetarandom} can be seen as a problem of optimization: indeed, let  $\mathcal G$ be the set of  probability densities on the interval
$(0,b),$ and consider the functional $\Psi : \mathcal G \longrightarrow \mathbb{R} ^ +$ defined, for any $g \in \mathcal G$, by
\begin{equation}
\Psi (g) = \left ( q- \int _ 0 ^b g(x) \pi_0(x) dx \right ) ^2 .
\end{equation}
Then, a solution $g$ to the IFPP problem, is characterized by
\begin{equation}
g = \arg \min _ { g \in \mathcal G } \Psi (g) .
\end{equation}
Of course, if there exists more than one  {density } $g \in \mathcal G$ at which $\Psi (g)$ attains the minimum, the solution of the IFPP problem \eqref{IFPPetarandom} is not unique; this also follows
by the fact that (see \eqref{maineq2}) the knowledge of $q$ provides only the mean of $\pi_0 (\eta),$ not its probability distribution.
\par\noindent
If one is looking for uniqueness of the solution to the IFPP problem \eqref{IFPPetarandom}, one must introduce constraints on the set $\mathcal G$ of probability densities of $\eta $ on the interval $(0,b)$
(see Example 2.1(iv), and \cite{abundo:CC22} and  {\cite{lefeb:22}}
in the case without resetting).
}
\end{Remark}

\begin{Remark}
{\it
{A solution to the IFPP problem may not even exist.
For instance, let us consider the case when $b=1$ and
$\pi_0(x)=1-x$ (see next Example 2.3 (i)), namely from \eqref{maineq2} $E(\eta) = 1-q,$  and let us
look for a solution $g$ in the class of unimodal densities in $(0,1)$ symmetric with respect to the middle point $1/2,$ that is, $E(\eta) = 1/2;$
if $q \neq 1/2, $  the equality $E(\eta) = 1-q,$ or equivalently Eq. \eqref{maineq2}, cannot be satisfied, and therefore the solution to the IFPP problem does not exist, under the above constraint on the density of $\eta .$
Of course, the same argument can be repeated for the density $h$ of the reset position $x_R,$ in the case when $x_R$ is not fixed, but random.}
}
\end{Remark}
{ As already mentioned, a diffusion process with resetting $\mathcal X(t)$ is able to modeling the time evolution of various systems in physics, biology, queueing theory, and other applied fields.
The  corresponding physical interpretation of the IFPP problem for such a process, from the experimental point of view, is as follows.
If the starting position $\eta=x \in (0,b)$ is fixed, but unknown, then
one can perform a large number of copies of observations of the physical process (theoretically described by $\mathcal X(t)$), in a sufficiently large time period $[0,T],$ and record the percentage of times
it first exits the strip $(0,b)$ from the lower side; thus, one obtains an empirical estimate $\hat \pi_0$ of the theoretical value $\pi_0.$ This result depends only on the intrinsic stochastic nature of the system; from $\hat \pi_0$ one has to find $x$ (in the theoretical description this is easy, because
Eq. (2.4) becomes  $q= \pi_0(x),$ being $g(y)= \delta (y-x),$ which is easily solved, since $\pi_0(x)$ is a monotone decreasing function). } \par\noindent
{ If $\eta$ is random, in addition to its intrinsic stochastic nature, the system is also subject to randomness due to the indeterminacy of the starting point; then,
by performing, as above, a large number of copies of observations of the process that starts its evolution from  $\eta ,$ one gets again an empirical estimate $\hat \pi_0$ of $q= \pi_0.$
Now, the goal is to retrieve the unknown probability distribution of the starting initial position $\eta.$
Actually, there are several scenarios in which this study is interesting; for instance, in biology in the context of stochastic models for neural activity with resetting.
The same considerations hold when one is looking for the probability distribution of $x_R$ (when the starting position $\eta=x$ and the reset rate $r$  are fixed, while the reset position $x_R$ is random). }

\bigskip

We start considering the case when the underlying diffusion $X(t)$ is BM  with drift $\mu,$ that is
$X(t)= x + \mu t + {W_t} \ (x = X(0)),$ then from \eqref{eqforpi0} with $Lf = \frac 1 2 f'' + \mu f'$ (see \eqref{generator}), we get that $\pi_0(x)$ is the solution of the problem:
\begin{equation} \label{eqpi0fordriftedBM}
\begin{cases}
\frac 1 2 f''(x) + \mu f'(x) -r f(x) = -rf(x_R), \ x \in (0,b)  \\
f(0)=1, \ f(b)=0.
\end{cases}
\end{equation}
By standard methods, one gets that the general solution of the  ODE in  \eqref{eqpi0fordriftedBM} is
$\pi_0 (x) = c_1 e^{d_1x} + c_2 e^{d_2 x} + \pi_0 (x_R),$ where $d_1 = - \mu - \sqrt {\mu ^2 + 2r } <0, \ d_2 = - \mu + \sqrt {\mu ^2 + 2r } >0, \
c_1, \ c_2$ are
arbitrary constants with respect to $x$ (they depend on $b, \ \mu, \ r$ and $x_R ).$
By imposing the boundary conditions $\pi_0 (0)=1, \ \pi_0(b)=0,$ one finally obtains:
$$ \pi_0(x_R) = -c_1 e^{d_1 b}- c_2 e^{d_2 b}$$
and
\begin{equation} \label{pi0fordriftedBM}
\pi_0(x) =
c_1 \left (e^ {d_1x} - e^{d_1b} \right ) + c_2 \left (e^ {d_2x} - e^{d_2b} \right ),
\end{equation}
where:
\begin{equation} \label{cifordriftedBM}
\begin{cases}
c_1 =  \left [1-\exp (d_1b ) - \exp (x_R(d_1-d_2)) \left ( 1- \exp (d_2b ) \right ) \right ] ^{-1}, \\
c_2 =  -  c_1 \exp (x_R (d_1 -d_2)) .
\end{cases}
\end{equation}
In the special case $\mu =0, \ d_1$ becomes $- \sqrt {2r}, \ d_2$ becomes $ \sqrt {2r},$    and the constants $c_i$ become:
\begin{equation} \label{pi0fornodriftBM}
\begin{cases}
c_1' =  \left [ 1- \exp (-b\sqrt {2r}) - \exp (-2x_R \sqrt {2r}) \left ( 1- \exp (b \sqrt {2r}) \right ) \right ] ^{-1}, \\
c_2' =   - c_1' \exp (-2x_R \sqrt {2r}).
\end{cases}
\end{equation}
For the sake of simplicity, for fixed $b>0$ we drop the dependence of the constants $d_i$ from $\mu$ and $r,$ and of $c_i$ and $c'_i$ from $ r$ and $x_R.$
\bigskip

\begin{Remark}
For \ $r=0$ (that is, when no resetting occurs), one obtains again the well-known result for drifted Brownian motion (without resetting), namely:
\begin{equation} \label{pi0r=0}
\pi _0(x) =
\begin{cases}
\frac { e^{-2 \mu x }- e^{-2 \mu b }} {  1- e^{-2 \mu b }}, & \mu \neq 0 \\
1- \frac x b, & \mu =0 .
\end{cases}
\end{equation}
\end{Remark}

If $g(x)$ is a solution to the IFPP problem \eqref{IFPPetarandom}, from \eqref{maineq}
it turns out that
\begin{equation} \label{valueofq}
q= \int _0 ^b g(x) \pi_0(x) dx = E[ \pi_0 ( \eta )]  =
\end{equation}
$$
 = \begin{cases}
c_1 \left (E[e^ {d_1 \eta}] - e^{d_1b} \right ) + c_2 \left (E[e^ {d_2 \eta}] - e^{d_2b} \right ), &  \mu \neq 0 \\
c_1' \left (E[e^ {- \eta  \sqrt {2r}}] - e^{ -b\sqrt {2r}} \right ) + c_2 ' \left (E[e^ {\eta \sqrt {2r}}] - e^{ b\sqrt {2r}} \right ), & \mu = 0,
\end{cases}
$$
where the constants $c_i, \ c_i'$ are defined in \eqref{cifordriftedBM}, \eqref{pi0fornodriftBM}, and $d_1 = - \mu - \sqrt {\mu ^2 + 2r }, \ d_2 = - \mu + \sqrt {\mu ^2 + 2r }.$
\par
We show below some explicit examples of solution to the IFPP problem  \eqref{IFPPetarandom}.
\bigskip

\noindent {\bf Example 2.1}
Let $\mathcal X(t)$ be the diffusion with resetting obtained from drifted BM $X(t)= \eta + \mu t + {W_t}.$ \bigskip

\noindent {\it (i)}
If
\begin{equation} \label{qex211}
q= \sum _{k=0} ^\infty \frac {c_1 d_1 ^k + c_2 d_2 ^k} {k!} \frac { B(\alpha +k, \beta )} {B( \alpha, \beta)} - c_1 e^{d_1b} - c_2 e^{d_2b};
\end{equation}
then, a solution $g$ to the IFPP problem \eqref{IFPPetarandom} for $\mathcal X(t),$ with $b=1,$ is the Beta density in $(0,1),$ namely
$g(x) = \frac { \Gamma (\alpha + \beta)} {\Gamma(\alpha) \Gamma ( \beta)} x^{\alpha -1} (1-x) ^{\beta -1} \mathbb{I}_{(0,1)}(x),$ where $\alpha, \ \beta >0.$
\par\noindent
To verify this by using Proposition \ref{mainpro} it is sufficient to substitute the various quantities into \eqref{valueofq}; it is convenient  to use  that, if $\eta$ has Beta density in $(0,1),$ one has (see e.g. \cite{gupta}):
\begin{equation} \label{BetaMGF}
 E[e^ {t \eta} ]= \sum _{k=0} ^\infty \frac {t^k} {k!} \frac { B(\alpha +k, \beta )} {B( \alpha, \beta)},
\end{equation}
where $$B( \alpha, \beta) = \frac {\Gamma (\alpha) \Gamma (\beta)} {\Gamma( \alpha + \beta)}.$$
\bigskip

\noindent If $\mu =0$ and $q$ is given by \eqref{qex211}, with $c_i$ replaced by $c'_i$ and with $d_1= - \sqrt {2r}, \ d_2=  \sqrt {2r},$ then
a solution $g$ to the IFPP problem for $\mathcal X(t)$  is again the Beta density in $(0,1).$ \par
{ If one takes $\alpha = \beta =1$ in \eqref{qex211}, that is
\begin{equation} \label{qEx12}
q= c_1 \left ( \frac 1 {d_1} \left (e^{d_1} -1 \right  ) -e^{d_1 b} \right )+ c_2 \left ( \frac 1 {d_2} \left (e^{d_2} -1 \right  ) -e^{d_2 b} \right ),
\end{equation}
then a solution $g$ to the IFPP problem \eqref{IFPPetarandom} for $\mathcal X(t)$ in the interval $(0,1)$ is the uniform density in $(0,1)$ (use \eqref{BetaMGF}). }
 \par\noindent
When e.g. $\mu =0$ and $r=1,$ the values of $q = q(x_R)$ obtained by \eqref{qEx12}  exhibit a decreasing behavior, as functions of $x_R;$ for instance: \par\noindent
$q(1/100) = 0.568; \ q(1/8)= 0.55; \ q(1/4)= 0.538; \ q(1/2)= 0.5; \ q(3/4)= 0.46;$
$q(0.9)= 0.441 .$
\par
{ If one takes $\alpha = \beta =2$ in \eqref{qex211}, that is
\begin{equation}
q= 6 \left \{ c_1 \left [ \frac {d_1 +2}  {d_1 ^3} + e^{d_1} \left ( - \frac 1 6 + \frac 1  {d_1^2} - \frac 2 { d_1^3} \right ) \right ]
+ c_1 \left [ \frac {d_2 +2}  {d_2 ^3} + e^{d_2} \left ( - \frac 1 6 + \frac 1  {d_2^2} - \frac 2 { d_2^3} \right ) \right ] \right \},
\end{equation} then a solution $g$ to the IFPP problem \eqref{IFPPetarandom} for $\mathcal X(t)$ in the interval $(0,1)$ is  $g(x)= 6 x (1-x) \mathbb{I}_{(0,1)}(x).$  }

\bigskip

\noindent {\it (ii)} Let $\mathcal X(t)$ be as in {\it (i)} and
let be
\begin{equation}
q= c_1 \left [ z(d_1) - e^{d_1} \right ] + c_2 \left [z(d_2) - e^{d_2} \right ],
\end{equation}
where $z(t):=\frac 4 {t^2} (e^{t/2} -1 )^2;$ then a solution $g$ to the IFPP problem \eqref{IFPPetarandom} for $\mathcal X(t)$ in the interval $(0,1)$
is the function $g$ defined by
\begin{equation} \label{funzgEx14}
g(x) =
\begin{cases}
4x, & x \in (0,1/2) \\
- 4x+4, & x \in [1/2,1) \\
0, & {\rm otherwise} .
\end{cases}
\end{equation}
In fact, by calculation, we get
$E(e^{t \eta}) = z(t),$ and the result soon follows, by inserting the various quantities into \eqref{valueofq}.
\bigskip

\noindent {\it (iii)}
Let $\mathcal X(t)$ be (undrifted) BM with resetting (that is the process in {\it (i)} with $\mu =0),$ and $b=1.$ If, for $\theta >0, \ \theta \neq \sqrt {2r} :$
\begin{equation}
q= c_1' \left ( \frac {\theta (e^\theta - e^{- \sqrt {2r}} ) } {(e^\theta -1) (\theta + \sqrt {2r} )}  - e^{- \sqrt {2r}} \right ) +
c_2' \left ( \frac {\theta (e^\theta - e^{ \sqrt {2r}} ) } {(e^\theta -1) (\theta - \sqrt {2r} )}  - e^{ \sqrt {2r}} \right ),
\end{equation}
then a solution $g$ to the IFPP problem \eqref{IFPPetarandom} for $\mathcal X(t)$ in the interval $(0,1)$ is the truncated exponential density in $(0,1)$ with parameter $\theta,$  that is
\begin{equation} \label{truncatedexp}
g(x) =
\begin{cases}
\frac {\theta e^{- \theta x } } {1- e^{- \theta} } , \ x \in (0,1), \\
0, \ \ \ \ \ \  {\rm otherwise}.
\end{cases}
\end{equation}
To verify this, it suffices to substitute the various quantities into \eqref{valueofq} and to use  that, if $\eta$ has density given by \eqref{truncatedexp},
one has:
\begin{equation}
E[e^{ - \lambda \eta } ] = \frac {\theta (e^\theta - e^{ - \lambda } ) } { (e^\theta -1) (\theta + \lambda )} , \ \lambda >0.
\end{equation}
\bigskip

\noindent {\it (iv)}
Let $\mathcal X(t)$ be drifted BM with resetting, as in {\it (i)}, and $b=1.$
Now, we look for solutions to the IFPP problem \eqref{IFPPetarandom} for $\mathcal X(t)$ of the form
\begin{equation} \label{formofg}
g(x) = a_1x + a_0, \ x \in (0,1),
\end{equation}
for suitable constants $a_0, \ a_1;$ notice that
$\int _0 ^1 g(x) dx =1$ implies  $a_1 /2 + a_0 =1.$  \par
By using Eq. \eqref{valueofq} and  taking into account that $a_1 /2 + a_0 =1,$  one obtains:
\begin{equation} \label{a0a1}
\begin{cases}
a_1=  \left \{ {q- \frac {c_1} {d_1}  (e^ {d_1} -1) - \frac {c_2} {d_2} (e^ {d_2} -1  ) +c_1 e^ {d_1} + c_2 e^ {d_2} }\right \} \times \\
\ \ \ \ \ \ \  \left \{ { \frac {c_1} {2d_1^2}  \left (e^ {d_1} (d_1 -2) + 2 +d_1 \right ) +
\frac {c_2} {2d_2^2} \left (e^ {d_2} (d_2 -2) + 2 +d_2 \right ) } \right \} ^{-1} \\
a_0 = 1- \frac {a_1} 2 .
\end{cases}
\end{equation}
Thus, for a given $q \in (0,1),$ one obtains that  $g(x) = a_1x + a_0$ with $a_0$ and $a_1$ given by \eqref{a0a1}, is the only solution
to the IFPP problem \eqref{IFPPetarandom} for $\mathcal X(t),$   with the constraint that solutions are searched in the set of linear functions in $(0,1)$ of the form \eqref{formofg}.
\bigskip

\noindent {\bf Example 2.2}
Let $\mathcal X(t)$ be the diffusion with resetting obtained from Geometric Brownian motion $X(t)$, which is driven by the SDE $(\theta, \ \sigma >0):$
\begin{equation}
 dX(t)= \theta X(t) dt + \sigma X(t) d{W_t}, \   X(0) = \eta ,
\end{equation}
whose solution is  $X(t)= \eta  e^{ \mu t +  \sigma {W_t}},$ with $\mu = \theta - \sigma ^2 /2.$ \par\noindent
Note that
$ Y(t) := \ln X(t) = \ln \eta + \mu t + \sigma {W_t} $ is drifted BM.
Then, for $\sigma =1$ the IFPP problem \eqref{IFPPetarandom} for $\mathcal X(t)$ in the interval $[1, b], \ b >1, $ becomes the IFPP problem for $\mathcal Y (t)$ in the interval $[0 , \ln b],$ where
$ \mathcal Y (t)$ is the diffusion with resetting obtained from $Y(t),$ with $y_R = \ln x_R .$ \par\noindent
Therefore, one easily obtains examples of solutions to the IFPP problem \eqref{IFPPetarandom} for $\mathcal X(t)$ in the interval $[1, b], \ b >1,$ by
using Examples 2.1 with suitable trivial modifications.
\bigskip

\noindent {\bf Example 2.3}
\par\noindent
{\it (i)}
Let $\mathcal X(t)$ be the diffusion with resetting obtained from the diffusion  $X(t),$
driven by the SDE
\begin{equation} \label{SDEex3}
dX(t)= \mu (X(t)) dt + \sigma (X(t)) d {W_t}, \ X(0) = \eta \in (0,b),
\end{equation}
where the drift $\mu (\cdot)$ and diffusion coefficient $\sigma (\cdot)$ are regular enough functions, such that there exists a unique strong solution of the
SDE.
As before, the probability $\pi _0(x)$ that $\mathcal X(t)$ first exits the interval $(0,b)$ from the left end, when starting from fixed $x \in (0,b),$ is the solution of
the differential equation with boundary conditions (see \eqref{eqforpi0}):
\begin{equation} \label{ex3}
\begin{cases}
\frac 1 2 \sigma ^2 (x) f'' (x) + \mu(x) f'(x) + r (f(x_R) - f(x) ) =0  \\
f(0)=1, \ f(b)=0 .
\end{cases}
\end{equation}
Now, we search for drift $\mu (x)$ and diffusion coefficient $\sigma (x),$ in such a way that $\pi _0(x)$ turns out to be equal to $1 - \frac x b,$ that is the probability that BM (without resetting) first exits the interval $(0,b)$ from the left end, when starting from $x \in (0,b);$ by substituting $f(x)= 1 - \frac x b$ in \eqref{ex3}, one gets
that it must be
$\mu(x) = r(x - x_R),$ whereas $\sigma (x)$ can be any, regular enough, diffusion coefficient.
\par\noindent
Therefore, we have obtained: \bigskip

let $\mathcal X(t)$ be the diffusion with resetting obtained from the diffusion  $X(t),$
driven by the SDE
\begin{equation}
dX(t)= r(X(t)-x_R) dt + \sigma (X(t)) d {W_t}, \ X(0) = \eta \in (0,b),
\end{equation}
where $\sigma (x)$ is any diffusion coefficient which satisfies the conditions for uniqueness of the solution of the SDE, and let be $q= \frac \beta {\alpha + \beta}, \ \alpha, \ \beta >0.$
Then, a solution to the IFPP problem \eqref{IFPPetarandom} for $\mathcal X(t)$ in the interval $(0,b)$ is the modified Beta density in $(0,b),$ namely
$$g(x)= \frac 1 {b^{\alpha + \beta -1}} \ \frac {x^ {\alpha -1} (b-x) ^ {\beta -1}} {B( \alpha, \beta)} \mathbb{ I}_{(0,b)}(x).$$
\par\noindent
This is soon verified by checking that \eqref{maineq} holds, being $g(x)$ the function above and $\pi_0(x) = 1 - x/b, \ x \in (0,b),$ or alternatively by using \eqref{maineq2}
and checking that it results $q= 1 - \frac 1 b E[ \eta] $ (use the expression of  the mean value of
the modified Beta random variable  $\eta$ in $(0,b),$  that is $E[ \eta] = \frac {b \alpha} {\alpha + \beta}).$
\bigskip

\noindent {\it (ii)}
Let us consider again the diffusion with resetting  $\mathcal X(t)$ obtained from the diffusion  $X(t),$ driven by Eq. \eqref{SDEex3}, and let be  $b=1;$  we search now for drift $\mu (x)$ and diffusion coefficient $\sigma (x),$ in such a way that
$\pi _0(x) = 2-2^x, \ x \in (0,1).$
By substituting $f(x)= 2-2^x$ in \eqref{ex3}, one gets that it must be
$$- \frac {(\ln 2 )^2} 2  \ \sigma ^2 (x)  - \ln 2 \ \mu(x) +r = r 2^{x_R-x}.$$
If e.g. $\sigma ^2(x) = \sigma ^2 = const,$ one obtains that
$\mu(x) = A -B / 2^{x},$
where $A= \frac r {\ln 2} - \frac  {\ln 2} 2 \sigma ^2, \
 B = \frac r {\ln 2}\ 2^{x_R}.$
 \par\noindent
Therefore, we have obtained: \bigskip

let $\mathcal X(t)$ be the diffusion with resetting obtained from the diffusion  $X(t)$  with drift $\mu(x)= A -B / 2^{x}$ and
diffusion coefficient $\sigma (x)= \sigma,$
and let
$$ q= 2- \sum _{k=0} ^ \infty \frac {(\ln 2) ^k } {k!} \ \frac {B( \alpha +k, \beta)}{B(\alpha, \beta)}= 1- \sum _{k=1} ^ \infty \frac {(\ln 2) ^k } {k!} \ \frac {B( \alpha +k, \beta)}{B(\alpha, \beta)} .$$
Then, a solution to the IFPP problem \eqref{IFPPetarandom} for $\mathcal X(t)$ in the interval $(0,1)$ is the Beta density in $(0,1),$ with parameters $\alpha$ and $\beta >0.$
\par\noindent
This is soon verified by checking that \eqref{maineq} holds, being $g(x)$ the Beta density and $\pi_0(x) = 2-2^x, \ x \in (0,1)$  (use also \eqref{BetaMGF} with
$t$ replaced by $\ln 2 ).$
\bigskip

\noindent {\bf Example 2.4}
 Let us suppose that the underlying diffusion $X(t)$ is Ornstein-Uhlenbeck (OU) process, which is driven by the SDE
 $dX(t)= - \nu X(t)dt + \sigma d {W_t} , $ with $\nu, \ \sigma >0,$ and
 let $\mathcal X(t)$ be the diffusion with resetting, obtained from the diffusion  $X(t).$ \par\noindent
Take e.g. $\sigma =1;$ then, the probability $ \pi _0 (x) $ that $\mathcal X(t)$ first exits the interval $(0,b)$ from the left end, when starting from fixed $x \in (0,b),$
satisfies the differential equation of the second order with non constant coefficients
\begin{equation} \label{p0OU_b=1}
\frac { 1} 2 f'' (x) - \nu x f '(x) - r f(x)= - rf(x_R) := {\kappa},
\end{equation}
with conditions  $f(0)=1, \ f(b)= 0.$ \par\noindent
Let $\tau_{OU}(x)$ be the first-passage time (FPT) of OU process $X(t)$ through zero, when starting from  $x >0,$ and denote by $\psi(x) = E[e^{-r\tau_{OU}(x)}]$ the Laplace transform of
$\tau_{OU}(x),$ which, for fixed $r>0,$  is a solution to the differential equation $\frac 12 \psi ''(x) - \nu x \psi '(x) - r \psi (x) =0,$ and it is explicitly given by Eq. (3.49) of
\cite{abundo:OUarea}, in terms of parabolic cylinder functions.
Then, in order to solve \eqref{p0OU_b=1},  we search for a solution of the form $f(x)= Z(x) \psi(x),$ where $Z(x)$ is a function to be found.
\par\noindent
By substituting this $f(x)$ into \eqref{p0OU_b=1} and taking into account that $\frac 12 \psi '' - \nu x \psi ' - r \psi =0,$ we finally obtain:
$$ \frac 1 2 \psi (x) Z''(x) + Z'(x) (\psi'(x) - \nu x \psi (x)) = {\kappa },$$
which is an ODE that contains only the first and second derivatives of $Z(x).$
In principle, this ODE can be solved by quadratures, so the function $Z(x)$ can be found, and consequently also $ \pi _0 (x).$ \par\noindent
Thus, if e.g. $q= \int _0 ^1 \pi_0 (x) dx =  \int _0 ^1 Z(x) \psi (x) dx ,$ then a solution to the IFPP problem \eqref{IFPPetarandom} for OU process with resetting in the interval $(0,1)$ is the uniform density in $(0,1).$
\bigskip

Of course, a solution $g$ to the IFPP problem \eqref{IFPPetarandom} can also be a discrete density, as
in the following examples. \bigskip

\noindent {\bf Example 2.5}
Let $\mathcal X(t)$ be BM  with resetting, and let be $0 < x_1 < b.$ If
\begin{equation}
q= \frac {c_1'} 3 \left ( 1 + e^{ - x_1 \sqrt {2r} } -2 e ^{ - b \sqrt {2r}}  \right ) +
\frac {c_2'} 3  \left ( 1 + e^{  x_1 \sqrt {2r} } -2 e ^{  b \sqrt {2r}} \right ),
\end{equation}
with $c_i'$ given by \eqref{pi0fornodriftBM},
then a solution $g$ to the IFPP problem \eqref{IFPPetarandom} for $\mathcal X(t)$ is the discrete uniform density on the set $\{ 0, x_1, b \}.$
To verify this, it suffices to substitute the various quantities in the second equation of \eqref{valueofq}, by taking into account that the Laplace transform of  the r.v. $\eta$ uniformly distributed on the set $\{ 0, x_1, b \}$ is
 $E[e^{- \lambda \eta }] = \frac 1 3 (1+ e^{- \lambda x_1} + e^{- \lambda b}), \ \lambda >0 .$
\bigskip

\noindent{\bf Example 2.6}
Let $b= n \ge 1$ be an integer and let $\mathcal X(t)$ be BM  with resetting.  \par\noindent
If, for $p \in (0,1):$
\begin{equation}
q= c_1 ' \left [\left (1-p+pe^{- \sqrt {2r}} \ \right )^n - e^{-n \sqrt {2r}} \right ] +  c_2 ' \left [ \left (1-p+pe^{ \sqrt {2r}} \ \right )^n - e^{n \sqrt {2r}} \right ],
\end{equation}
with $c_i'$ given by \eqref{pi0fornodriftBM},
then a solution $g$ to the IFPP problem \eqref{IFPPetarandom} for $\mathcal X(t)$ in the interval $[0,n]$ is the Binomial density with parameters $n$ and $p.$
To obtain this, it suffices to substitute the various quantities in \eqref{valueofq}, by taking into account that,
if $\eta$ has binomial distribution  with parameters $n$ and $p$, then its Laplace transform
is $E[e^{- \lambda \eta }] = (1-p+p e^{- \lambda } )^n, \ \lambda >0 .$
\bigskip

\noindent {\bf Diffusion processes conjugated to BM.} \ We recall that a one-dimensional diffusion process $X(t)$ starting from $X(0) = x $ is said to be conjugated to Brownian motion, if there exists an increasing, differentiable function $v(x)$
with $v(0)=0$ such that $X(t) = v^{-1} ({W_t} + v ( x))$ (see e.g.  \cite{abundo:stapro12}).
\par\noindent
For instance: \par\noindent
$\bullet$ the Feller process $X(t)$ driven by the SDE: $$dX(t) = \frac 1 4 dt + \sqrt {X(t)} \ d {W_t}, \ X(0) = x \ge 0,$$ is conjugated to BM via the function $v(x) =2 \sqrt x,$
that is $X(t) = \frac 1 4 ({W_t} + 2 \sqrt x )^2;$ \par\noindent
$\bullet$ the Wright$\&$Fisher-like process $X(t)$ driven by the SDE:
$$dX(t)= \left (\frac 1 4 - \frac  1 2 X(t) \right ) dt + \sqrt {X(t)(1-X(t)} \ d {W_t}, \ X(0) \in (0,1),$$ is conjugated to BM via the function
 $v(x)= 2 \arcsin \sqrt x,$
that is:
\par\noindent
 $X(t) = \sin ^2 \left ( \frac 1 2 {W_t} + \arcsin \sqrt x \right ) .$ \bigskip

\noindent If the diffusion $X(t)$  is conjugated to BM  via the function $v,$ the  corresponding process $\mathcal X(t),$  which is obtained by resetting $X(t)$ to  $x_R \in (0,b)$  at the rate $r,$ is transformed  via the function $v(x)$ into  BM with resetting with  $v(x_R)$ in place of $x_R.$ Note that, for fixed $x \in (0,b), \pi_0 (x) = P(\mathcal X(\tau _ {0,b}(x))=0)$ is nothing but  $ \pi _0 ^B(v(x)) ,$ that is, the probability that BM with resetting starting from $v(x)$ first exits the interval
$(0, v(b))$ through the left end $0;$ so, by replacing $b$ with $v(b)$ and $x_R$ with $v(x_R),$ and taking $\mu =0$ in \eqref{pi0fordriftedBM} we obtain:
\begin{equation}
\pi_0(x) =
\bar c_1' \left (e^ {-v(x) \sqrt {2r}} - e^{ -v(b)\sqrt {2r}} \right ) + \bar c_2 ' \left (e^ {v(x) \sqrt {2r}} - e^{ v(b)\sqrt {2r}} \right ),
\end{equation}
where (see \eqref{pi0fornodriftBM}):
\begin{equation} \label{constantsbarci'}
\bar c_1' = \frac 1 {1-e^{-v(b)\sqrt {2r}} - e^{-2v(x_R) \sqrt {2r}} \left ( 1- e^{v(b) \sqrt {2r}} \right ) }, \ \bar c_2' = - \bar c_1 ' e^{-2 v(x_R) \sqrt {2r}}.
\end{equation}
Of course, if the starting point $X(0)$ is not deterministic, but it is a random variable $\eta \in (0,b),$
one has $v(X(t)) = {W_t} +  \tilde \eta,$ where $\tilde \eta = v( \eta) \in (0, v(b))$ is the random starting point of the corresponding BM.
Then, if $\tilde g$ is the density of $\tilde \eta,$ one gets that  $\eta$ has density $g(x) = \tilde g(v(x)) v'(x), \ x \in (0, b).$
\bigskip

\noindent {\bf Example 2.7}
Let us suppose that the underlying process $X(t)$ is conjugated to BM via a function $v(x)$ with $v(b) =1,$ and let
\begin{equation}
q= \sum _{k=0} ^\infty \frac {\bar c_1' (- \sqrt {2r}) ^k + \bar c_2' ( \sqrt {2r}) ^k} {k!} \frac { B(\alpha +k, \beta )} {B( \alpha, \beta)} - \bar c_1' e^{-b\sqrt {2r}} - \bar c_2' e^{b \sqrt{2r}},
\end{equation}
where $\bar c_i '$ are given by \eqref{constantsbarci'};
then, a solution $g$ to the IFPP problem \eqref{IFPPetarandom} for $\mathcal X(t)$ in the interval $(0,1)$ is the density obtained by transformation via $v^ {-1}$ of the Beta density in $(0,1),$
namely
$$g(x) = \frac { \Gamma (\alpha + \beta)} {\Gamma(\alpha) \Gamma ( \beta)} v(x)^{\alpha -1} (1-v(x)) ^{\beta -1} v'(x) \mathbb{I}_{(0,1)}(x)$$
(see Example 2.1(i) with $\mu =0)$.
In particular, if
$$ q=   \bar c_1' \left ( - \frac 1 {\sqrt {2r}} \left (e^{- \sqrt {2r}} -1 \right ) - e^{- \sqrt {2r}} \right ) + \bar c_2 ' \left (  \frac 1 {\sqrt {2r}} \left (e^{ \sqrt {2r}} -1 \right ) - e^{ \sqrt {2r}} \right ),$$
then, a solution $g$ to the IFPP problem \eqref{IFPPetarandom} for $\mathcal X(t)$ in the interval $(0, 1)$ is obtained by transformation via $v^ {-1}$ of the uniform density in $(0,1),$ i.e.
$g(x)=  v'(x) \mathbb{I}_{(0,1)}(x)$ (cf. \eqref{qEx12}).
\section{Case I: Inverse first-passage time (IFPT) problems for random initial position $\eta$ and fixed $r$ and $x_R$}
In this section we suppose that $\mathcal X(t)$ is  drifted BM with resetting (with drift $\mu).$
\subsection{The IFPT problem}
We recall the terms of the IFPT problem \eqref{IFPTetarandom}.
Let $\mathcal X(t)$ be drifted BM with resetting,
starting from a random initial position $\eta \in (0, + \infty),$
which is supposed to be independent of $\mathcal X(t)$.
Let $\tau(x)$ be the first-passage time (FPT) of $\mathcal X(t)$ through zero, under the condition that $\eta  =x >0,$ namely
\begin{equation}
\tau (x) = \inf \{ t>0: \mathcal X(t) =0  | \eta = x \},
\end{equation}
and let
$\tau$ the (unconditional) FPT of $\mathcal X(t)$ through zero, i.e.
\begin{equation}
\tau = \inf \{ t>0: \mathcal X(t) =0 \}.
\end{equation}
For a given  distribution function $F (t)$ on the positive real axis, or equivalently for a given density $f(t) = F'(t), \ t >0,$ the IFPT problem
  \eqref{IFPTetarandom} consists in finding the density $g$ of the
random initial position $\eta \in (0, + \infty),$  if it exists, such that $P(\tau \le t ) = F(t), \ t >0.$ The function $g$ is called a solution to the IFPT problem
\eqref{IFPTetarandom} for $\mathcal X(t).$ \par\noindent
We recall from \cite{abundo:FPA2023}, Eq. (4.1), that for fixed $x \in (0, + \infty)$ the Laplace transform (LT) of $\tau (x)$ is, for $\lambda >0:$
\begin{equation} \label{LTtaux}
\widehat f _ { \tau(x)} (\lambda |x) = E \left [ e^{- \lambda \tau (x)} \right ]= e^{ -x \left ( \mu + \sqrt {\mu ^2 + 2(\lambda + r)} \ \right )}
\end{equation}
$$ + \frac r {\lambda +r e^{- x_R \left ( \mu + \sqrt {\mu ^2 + 2(\lambda + r)} \ \right )}}  \left ( 1- e^{-x \left ( \mu + \sqrt {\mu ^2 + 2(\lambda + r)} \ \right )} \right )
 e^{- x_R \left ( \mu + \sqrt {\mu ^2 + 2(\lambda + r)} \ \right )}
$$
$$
= e^{ -x \left ( \mu + \sqrt {\mu ^2 + 2(\lambda + r)} \ \right )} +
C(\lambda, \mu, x_R) \left ( 1- e^{-x \left ( \mu + \sqrt {\mu ^2 + 2(\lambda + r)} \ \right )} \right ),
$$
where, for fixed $r>0:$
\begin{equation} \label{ClammuxRbis}
C(\lambda, \mu, x_R) =  \frac {r e^{- x_R \left ( \mu + \sqrt {\mu ^2 + 2(\lambda + r)} \ \right )}} {\lambda +r e^{- x_R \left ( \mu + \sqrt {\mu ^2 + 2(\lambda + r)} \ \right )}}.
\end{equation}
Then, if $g(x)$ denotes the density of the random initial position $\eta ,$ the LT of $\tau$ turns out to be:
\begin{equation}
\widehat f _ { \tau} (\lambda ) = E \left [ e^{- \lambda \tau } \right ]= \int _0 ^{ + \infty } g(x) E \left [e^{- \lambda \tau (x)} \right ] dx
\end{equation}
$$
= \int _0 ^{+ \infty } g(x) e^{-x \left ( \mu + \sqrt {\mu ^2 + 2(\lambda + r)} \ \right )} \left (1- C(\lambda, \mu, x_R) \right ) dx + C(\lambda, \mu, x_R)
$$
$$
=
\left (1- C(\lambda, \mu, x_R) \right ) \widehat g \left ( \mu + \sqrt {\mu ^2 + 2 (\lambda +r)} \right ) + C(\lambda, \mu, x_R),
$$
where $\widehat g (\theta)= \int _0 ^{+ \infty } e^{- \theta x} g(x) dx, \ \theta >0, $ is the Laplace transform of $g(x).$ \par\noindent
Thus, we have obtained the following:

\begin{Proposition}
Let $\mathcal X(t)$ be drifted BM with resetting, starting from the random initial position $\eta \in (0, + \infty),$ which is supposed to be independent of $\mathcal X(t),$ and let $f(t), \ t >0,$ be a given  density. Then, if there
exists a  solution $g$  to the IFPT problem \eqref{IFPTetarandom} for $\mathcal X(t),$ the following equality holds between the LTs  $\widehat f$ of $f$ and the LT $\widehat g$ of $g:$
\begin{equation} \label{LTtau}
\widehat f  (\lambda ) = \left (1- C(\lambda, \mu, x_R) \right ) \widehat g \left ( \mu + \sqrt {\mu ^2 + 2 (\lambda +r)} \right ) + C(\lambda, \mu, x_R), \ \lambda >0 ,
\end{equation}
where
\begin{equation} \label{ClambdamuxR}
C(\lambda, \mu, x_R) =  \frac {r e^{- x_R \left ( \mu + \sqrt {\mu ^2 + 2(\lambda + r)} \ \right )}} {\lambda +r e^{- x_R \left ( \mu + \sqrt {\mu ^2 + 2(\lambda + r)} \ \right )}}.
\end{equation}
\end{Proposition}
\hfill  $\Box$
\begin{Remark}
From \eqref{LTtau}, one gets for $|\theta - \mu | \neq  \sqrt {\mu ^2 + 2r}:$
\begin{equation}
\widehat g ( \theta) =
\frac 2 {\theta ^2 -2r - 2 \theta \mu }
\left [  \left ( \frac {\theta ^2} 2 -r - \theta \mu   +r e^ {- \theta x_R }  \right ) \widehat f \left (\frac {\theta ^2} 2 -r - \theta \mu \right )
-r e^ {- \theta x_R } \right ] .
\end{equation}
Therefore, the LT of $\eta$ exists in the open interval, containing the origin  \par\noindent
$ \left ( \mu - \sqrt {\mu ^2 + 2r}, \ \mu + \sqrt {\mu ^2 + 2r} \right ),$ so it is infinite times differentiable
at $\theta =0,$ and there exist the moments of $\eta$ of any order. If $\widehat g( \theta)$ is an analytic function,
 it uniquely identifies the distribution of $\eta;$ in this case, unlike the case of the IFPP problem \eqref{IFPPetarandom}, if a solution to the IFPT problem \eqref{IFPTetarandom} exists, it is unique.

\end{Remark}

We show below some examples of solution to the IFPT problem \eqref{IFPTetarandom}; for the sake of simplicity, we limit ourselves to the case when $\mathcal X(t)$ is (undrifted) BM with resetting $(\mu =0).$
\bigskip

\noindent {\bf Example 3.1} Let $\mathcal X(t)$ be BM with resetting, starting from the random initial position $\eta \in (0, + \infty).$ \par\noindent
For $\nu >0,$ let us suppose that the  density $f$  has LT:
\begin{equation} \label{exampleLTFPTdensity}
\widehat f  (\lambda ) = \frac 1 {\lambda +r e^{- x_R \sqrt {2(\lambda +r)}}} \left ( \frac {\lambda \nu} {\nu + \sqrt {2(\lambda +r)}} + r e^{-x_R \sqrt {2(\lambda +r)}}  \right )
 , \ \lambda >0 .
\end{equation}
Then, the solution $g$ to the IFPT problem \eqref{IFPTetarandom} for the  density corresponding to the  LT \eqref{exampleLTFPTdensity} is the exponential density with parameter $\nu,$ that is $g(x)= \nu e^{- \nu x}, \ x >0 .$
To verify this it is sufficient to recall that the LT of this density $g$ is $\widehat g (\theta)= \nu / (\nu + \theta), \ \theta >0,$ and to use \eqref{LTtau}.
The density $f$ corresponding to the LT \eqref{exampleLTFPTdensity} cannot be obtained in closed form.
However, we can get some qualitative characteristics of the distribution having density $f;$ in fact, from $\widehat f  (\lambda )$ we easily obtain all the moments $m_k$ and the central ones $\mu_k \ , k=1, 2 , \dots  $ of the distribution. If e.g. $\nu =r = x_R =1,$ the mean of the distribution turns out to be $2.41, $ and the first three central moments are:
$$  \mu_2 = 9.61, \ \mu _3 = - 66.07 , \ \mu _4= 485.81 , $$
from which skewness  $\gamma _1 = \frac {\mu _3 } {\mu_2 ^{3/2}}  = -2.217,$ and excess kurtosis coefficient $\gamma _2 = \frac {\mu _4} { \mu _2 ^2 } -3=
2.26$ follow. Since skewness $\gamma _1 $ is negative,
the tail of the distribution is on the left side;
moreover, the density $f(t)$ tends to zero, as $ t \rightarrow + \infty,$ more slowly than the normal density does,  because excess kurtosis $\gamma _2 $ is positive. \par
 \par
As an additional example, if for $\alpha, \ \nu > 0, $ the LT of the  density $f$ is:
\begin{equation} \label{example2LTFPTdensity}
\widehat f  (\lambda ) = \frac 1 {\lambda +r e^{- x_R \sqrt {2(\lambda +r)}}} \left [ \frac {\lambda \nu ^\alpha } { \left (\nu + \sqrt {2(\lambda +r)} \ \right )^\alpha } + r e^{-x_R \sqrt {2(\lambda +r)}}  \right ], \ \lambda >0
\end{equation}
then, the solution $g$ to the IFPT problem \eqref{IFPTetarandom} for the density corresponding to the LT \eqref{example2LTFPTdensity} is the Gamma density with parameters $\alpha$ and $\nu.$
This soon follows by using that the LT of this Gamma density is $\widehat g (\theta)= [\nu / (\nu + \theta)] ^ \alpha , \ \theta >0.$
The previous example is a special case, when $\alpha =1.$
\bigskip

The following examples concern discrete densities $g,$ as  solutions to the IFPT problem \eqref{IFPTetarandom} for BM with resetting. \bigskip

\noindent{\bf Example 3.2}
Let $\mathcal X(t)$ be BM with resetting, starting from the random initial position   $\eta ,$ and assume that $\eta$ is an non-negative integer.
For $p \in (0,1),$ let us suppose that the density $f$  has LT:
\begin{equation} \label{exampleLTFPTdensity2}
\widehat f  (\lambda ) = \left (1- C(\lambda, 0, x_R) \right ) \frac p {1- (1-p) e^{ - \sqrt {2(\lambda +r) } }  } + C(\lambda, 0, x_R)
 , \ \lambda >0 ,
\end{equation}
where $C(\lambda, \mu, x_R)$ is given by \eqref{ClambdamuxR}. \par\noindent
Then, the solution $g$ to the IFPT problem \eqref{IFPTetarandom} for the  density corresponding to the LT  \eqref{exampleLTFPTdensity2} is the Geometric density with parameter $p,$
that is, $g(k)= P(\eta =k )= p(1-p) ^k, \ k=0, 1, \dots .$
This soon follows by using that the LT of this Geometric density is $\widehat g (\lambda)=  \frac p {1-(1-p) e^{- \lambda }}, \ \lambda >0.$ \bigskip

\noindent{\bf Example 3.3} Let $\mathcal X(t)$ be BM with resetting; for $\nu >0,$ suppose that
the density $f$  has LT:
\begin{equation} \label{exampleLTFPTdensity3}
\widehat f  (\lambda ) = \left (1- C(\lambda, 0, x_R) \right ) \exp \left [ \nu \left ( e^{- \sqrt {2(\lambda +r ) } } -1 \right )  \right ]   +C(\lambda, 0, x_R)
 , \ \lambda >0 .
\end{equation}
Then, the solution $g$ to the IFPT problem \eqref{IFPTetarandom} for the density having LT given by \eqref{exampleLTFPTdensity3} is the Poisson density with parameter $\nu,$
that is, $g(k)= P(\eta =k )= e^{-\nu } \frac { \nu ^k} {k! },  \ k=0, 1, \dots .$
It suffices to use that the LT of this Poisson density is $\widehat g (\lambda)=  \exp \{ \nu (e^{- \lambda } -1 )  \} , \ \lambda >0$ \bigskip

\subsection{The Inverse Mean-FPT (IMFPT) problem}
Let $\mathcal X(t)$ be drifted BM with resetting (with drift $\mu ),$ starting from the random initial position $\mathcal X(0) = \eta >0,$
which is independent from $\mathcal X(t),$ where  $r$ and $x_R >0$ are fixed.
If $ \tau  = \inf \{ t>0: \mathcal X(t) =0 \}$ is the FPT of $\mathcal X(t)$ through zero, one has (see \cite{abundo:FPA2023}, Eq. (4.3)):
\begin{equation} \label{mean2}
E [ \tau (x) ] := E [ \tau | \eta =x]= \frac 1 r \left(1- e^{-x \left (\mu +  \sqrt {\mu ^2 + 2r} \ \right )} \right ) e^{ x_R \left (\mu +  \sqrt {\mu ^2 + 2r} \ \right ) }.
\end{equation}
Thus, if $\eta$ has density $g,$  the (unconditional) mean of the FPT  becomes:
\begin{equation} \label{Etauetarandom}
E[ \tau ] = \frac  1 r e^{ x_R \left (\mu +  \sqrt {\mu ^2 + 2r} \ \right ) } \int _0 ^{+ \infty} g(x)  \left [ 1- e^{-x  \left (\mu +  \sqrt {\mu ^2 + 2r} \  \right ) } \right ] dx \end{equation}
$$
= \frac  1 r e^{ x_R \left (\mu +  \sqrt {\mu ^2 + 2r} \ \right ) } \left [1- \widehat g \left ( \mu + \sqrt {\mu ^2 +2r} \ \right ) \right ] ,
$$
where $\widehat g $ denotes the Laplace transform of $g.$
\par
Now, we consider the following kind of inverse problem for the mean of the FPT (IMFPT problem):
\begin{equation} \label{IMFPTcaseI}
{given \ m >0, \ find \ the \ density \ g \ of \ \eta \ (if \ it \ exists),\ so \ that \ E[\tau] =m}
\end{equation}
(note that $E[\tau]$ is indeed a function of $\mu, \ x_R $ and $r$).  The {density} $g$ is called a solution to the inverse problem
\eqref{IMFPTcaseI}; also now, uniqueness of the solution is not guaranteed.\bigskip

\noindent {\bf Example 3.4}
Let $\mathcal X(t)$ be drifted BM with resetting, starting from the random initial position  $\eta >0,$  and for $a, \theta >0$  let be
$$ m = \frac 1 r e^{ x_R \left (\mu +  \sqrt {\mu ^2 + 2r} \ \right ) } \left [ 1-  \left ( \frac \theta {\theta + \mu + \sqrt {\mu ^2 + 2r}}\right ) ^a \right ].$$
Then, a solution $g$ to the IMFPT problem  \eqref{IMFPTcaseI} is the Gamma density with parameters $a$ and $\theta .$ In particular, if
$$ m = \frac 1 r e^{ x_R \left (\mu +  \sqrt {\mu ^2 + 2r} \ \right ) } \left [ 1-   \frac \theta {\theta + \mu + \sqrt {\mu ^2 + 2r}} \right ],$$
a solution $g$ is the exponential density with parameter $\theta.$
\par\noindent
To verify this, it suffices to use the expression of the LT of a Gamma density and to insert the various quantities into \eqref{Etauetarandom}.
\bigskip

\subsection{The Inverse Mean-FET (IMFET) problem}
Let $\mathcal X(t)$ be drifted BM with resetting (with drift $\mu ),$ starting from the random initial position $\eta \in (0,b),$ which is independent from
$\mathcal X(t)$, where  $r$ and $x_R \in (0,b)$ are fixed.
If
$\tau _{0,b} = \inf \{ t>0: \mathcal X(t) \notin (0,b) \}$ denotes the first-exit time (FET) of $\mathcal X(t)$ from $(0,b),$ then
$E[\tau _{0,b}(x)] := E[\tau _{0,b} | \eta =x ],$ as a function of $x,$ is the solution of the differential equation with boundary conditions (see e.g. \cite{abundo:pms00}):
\begin{equation}
\begin{cases}
\frac 1 2 f''(x) + \mu f'(x) - r f(x) + r f(x_R) = -1 \\
f(0) = f(b) =0.
\end{cases}
\end{equation}
By solving, one finds:
\begin{equation}
E[\tau _{0,b}(x)] =  C_1 (e^{d_1 x} -1) + C_2 (e^{d_2 x} -1),
\end{equation}
where $d_1 = - \mu - \sqrt {\mu ^2 + 2r}, \ d_2 = - \mu + \sqrt {\mu ^2 + 2r} $ and
\begin{equation} \label{costantsC}
\begin{cases}
C_1 =  \frac 1 r e^{- d_2 x_R } (1-e^{d_2 b}) \left [ 1- e^ {d_1b} - e^{x_R (d_1 - d_2)} (1- e^{d_2b} ) \right ] ^{-1}
 \\
C_2 = - \frac 1 r e^{- d_2 x_R } (1-e^{d_1 b})  \left [ 1- e^ {d_1b} - e^{x_R (d_1 - d_2)} (1- e^{d_2b} ) \right ] ^{-1}
\end{cases}
\end{equation}
(note that the constants $C_i$ depend on $b, \ \mu, \ r$ and $x_R).$ \par
If the density of   $\eta$ is  $g(x), \ x \in (0,b),$ then the mean FET becomes
\begin{equation} \label{meanFETetarandom}
E[\tau _{0,b}] = \int _0 ^b g(x) E[\tau _{0, b} (x)] dx =   C_1 (E[e^{d_1 \eta}] -1) + C_2 (E[e^{d_2 \eta}] -1).
\end{equation}
Now, one can consider the following kind of inverse problem for the mean of the FET (IMFET problem):
\begin{equation} \label{IMFETcaseI}
{given \ m >0, \ find \ the \ density \ g \ of \ \eta \ (if \ it \ exists),\ so \ that \ E[\tau _{0,b}] =m.}
\end{equation}
The {density $g$} is called a solution to the inverse problem \eqref{IMFETcaseI}; in fact, also in this case the uniqueness of solution is not guaranteed.
\bigskip

\noindent {\bf Example 3.5}
If $b=1$ and, for $\alpha, \beta >0:$
$$ m = C_1 \left [  \sum _ {k=1} ^\infty \frac { d_1 ^k} {k! }  \ \frac {B(\alpha + k, \beta) } { B(\alpha, \beta)}  \right ] +
C_2 \left [  \sum _ {k=1} ^\infty \frac { d_2 ^k} {k! }  \ \frac {B(\alpha + k, \beta) } { B(\alpha, \beta)}  \right ],$$
then, a solution $g$ to the inverse problem \eqref{IMFETcaseI} is the Beta density in $(0,1).$ \par\noindent
This soon follows by calculating $C_1 (E[e^{d_1 \eta}] -1) + C_2 (E[e^{d_2 \eta}] -1),$ and making use of  \eqref{meanFETetarandom} and \eqref{BetaMGF}.\par\noindent
If $\alpha = \beta =1,$ a solution is the uniform density in $(0,1).$
\bigskip

The following is  an example in which a solution to the IMFET problem \eqref{IMFETcaseI} is a discrete density. \bigskip

\noindent {\bf Example 3.6}
For $b$ equal to an integer $n,$ let  $\mathcal X(t)$ BM with resetting, and
$$m= C_1 \left [ (1-p+pe^{ - \sqrt {2r}})^n -1 \right ] + C_2 \left [ (1-p+pe^{  \sqrt {2r}})^n -1 \right ] ,$$
for some $p \in (0,1),$
where $C_1, C_2$ are  obtained from  \eqref{costantsC}, by taking   $b=n$ and  $d_1 = - \sqrt {2r}, \ d_2 = \sqrt {2r}$ (because $\mu =0) .$ \par\noindent
Then, a solution $g$ to the IMFET problem \eqref{IMFETcaseI} is the Binomial density with parameters $n$ and $p.$
This easily follows by using that, if $\eta $ has such Binomial distribution, its  LT is $E[ e^{-\lambda \eta}] = (1-p+pe^{-\lambda} )^n, \ \lambda >0.$
\bigskip

In the next Section, we briefly address IFPP, IFPT, IMFPT, and IMFET problems considered so far for  the process $\mathcal X(t)$ with resetting, in the case II, namely when the resetting rate $r>0$ and the initial position $\eta = \mathcal X(0)$ are fixed, whereas  the reset position $x_R$ is random.

\section {Case II: Inverse problems when the reset position $x_R$ is random, whereas $r$ and $\mathcal X(0)$ are fixed}
In this section, $\mathcal X(t)$ is Wiener process with resetting.
\subsection{The IFPP problem}
Let $\mathcal X(t)$ be drifted BM with resetting; $ r >0$ and $\eta = \mathcal X(0) =x \in (0,b)$ are  fixed, whereas   $x_R \in (0,b)$ is random and independent of $\mathcal X(t);$
the corresponding IFPP problem consists in the following:
\begin{equation} \label{IFPPcaseII}
{given \ q \in (0,1), \ find \ the \ density \ h \ of \ x_R \ (if \ it \ exists),\ so \ that \ \tilde \pi _0 =q,}
\end{equation}
where:
$$\tau _ {0,b} (u) = \inf \{ t >0: \mathcal X(t) \notin (0,b) | x_R =u \}$$
is the FET of $\mathcal X(t)$ from $(0,b),$ conditional to $x_R=u \in (0,b),$
$$ \tilde \pi _0(u) = P(\mathcal X(\tau _ {0,b} (u)) =0 |x_R =u)$$
is the probability that $\mathcal X(t)$ first exits the interval $(0,b)$ from the left end, under the condition that $x_R=u,$
 and
 $$\tilde \pi _0 = \int _0 ^ b h(u) \tilde \pi _0(u) du $$
 is the (unconditional) probability that $\mathcal X(t)$ first exits the interval $(0,b)$ from the left end. \par\noindent
 The {density $h$} is called a solution to the IFPP problem \eqref{IFPPcaseII}. \par\noindent
From \eqref{pi0fordriftedBM}, one gets that a solution $h$ to the IFPP problem \eqref{IFPPcaseII} satisfies the integral equation:
\begin{equation} \label{integraleqforqtwo}
q= \int _0 ^b h(x_R) \left [ c_1 \left ( e^{d_1 x} - e^{d_1b} \right ) + c_2 \left (e^{d_2 x} - e^{d_2b} \right ) \right ] dx_R,
\end{equation}
where $d_i$ are as in Section 2 and $c_i$ are given by \eqref{cifordriftedBM} (note that now, unlike the case of Section 2,  $c_i$ have to be considered as functions of $r$ and $x_R,$ being $x$ fixed). \bigskip

\noindent
{\bf Example 4.1}
Let be $\mathcal X(t)$ (undrifted) BM with resetting, and $b=1;$ moreover, let be:
\begin{equation}
A= e^{-x \sqrt {2r}} - e^{- \sqrt {2r}}, \ B= e^{x \sqrt {2r}} - e^{\sqrt {2r}}, \ C= 1- e^{- \sqrt {2r}}, \ D= 1- e^{ \sqrt {2r}}
\end{equation}
and
\begin{equation}
\alpha = \frac A C + \frac {DB } {C^2 }, \
\beta = \frac B C ,   \
\gamma =  \frac {AD } C + \frac {B D^2 }  {C^2 }  \ .
\end{equation}
If
\begin{equation}
q= \frac 1 { 2 \sqrt {2r} } \Big [2 \alpha \sqrt {2r} - \beta e^{ - 2 \sqrt {2r}} -  \frac \gamma D \ \ln \left (C - D e ^{2 \sqrt {2r}} \right )
 -  \beta - \frac \gamma D \ \ln \left (C - D  \right ) \Big ],
\end{equation}
then a solution $h$ to the IFPP problem \eqref{IFPPcaseII} is the uniform density in $(0,1).$ To verify this, it suffices to use \eqref{pi0fornodriftBM} and to substitute the various quantities into \eqref{integraleqforqtwo}.

\subsection{The IFPT problem}
Let $\mathcal X(t)$ be drifted BM with resetting; we suppose that $r >0$ and $\eta = \mathcal X (0) = x \in (0, + \infty)$ are fixed, whereas $x_R \in (0, + \infty)$ is random and independent of $\mathcal X(t) .$
In this situation, if $\tau$ is the FPT of $\mathcal X(t)$ through zero, the corresponding IFPT problem consists in the following:
\begin{equation} \label{IFPTxRrandom}
\begin{aligned}
given \ a \ distribution \ function \ F(t) \ on \ the \ positive \ real \ axis, \ find   \\
the \ density \ h \ of \ x_R >0 \ (if \ it \ exists), \ such \ that \ P( \tau \le t ) = F(t).
\end{aligned}
\end{equation}
From \eqref{LTtaux}, one gets that $h(u)$ is a solution to this IFPT problem, if  the LT of $f(t):= F'(t),$ i.e. $\widehat f (\lambda)= E[e^{ - \lambda \tau }],$ satisfies
\begin{equation} \label{LTtwo}
\widehat f (\lambda)  = e ^{ - x ( \mu + \sqrt {\mu ^2 + 2 ( \lambda +r) } \ )  }
\end{equation}
$$+ (1 - e ^{ - x ( \mu + \sqrt {\mu ^2 + 2 ( \lambda +r) } \ )  }  )
\int _ 0 ^{+ \infty} \frac { r h(u) e^ {- u ( \mu + \sqrt {\mu ^2 + 2 ( \lambda +r) } \ ) } }  { \lambda + r e^ {- u ( \mu + \sqrt {\mu ^2 + 2 ( \lambda +r) } \ ) } } du  , \ \lambda >0.
$$
We show below a simple example of solution to the IFPT problem \eqref{IFPTxRrandom}, in the case of undrifted BM with resetting $(\mu =0)$.
\bigskip

\noindent {\bf Example 4.2}
Let $\mathcal X(t)$ be  (undrifted) BM with resetting and
suppose that  the LT of the density $f$ is:
\begin{equation}
\widehat f (\lambda) =e ^{ - x \sqrt {2 ( \lambda +r) } } + (1 - e ^{ - x \sqrt {2 ( \lambda +r) } }  )
\frac 1 {x \sqrt {2(\lambda +r)}} \ln \left ( \frac {\lambda +r } {\lambda +r e^{-x \sqrt {2(\lambda +r)}} } \right )
, \ \lambda >0.
\end{equation}
Then, the solution $h$ to the IFPT problem \eqref{IFPTxRrandom} is the uniform density in $(0,x).$ To verify this, it is sufficient to substitute the various quantities into \eqref{LTtwo}.

\subsection{The Inverse Mean-FPT (IMFPT) problem}
Let $\mathcal X(t)$ be drifted BM with resetting (with drift $\mu);$ we suppose that $r >0$ and $\eta = \mathcal X(0) =x \in (0, + \infty)$ are  fixed, whereas  $x_R \in (0, + \infty)$ is random and independent of $\mathcal X(t).$
Let
$ \tau  = \inf \{ t>0: \mathcal X(t) =0 \}$ be the FPT of $\mathcal X(t)$ through zero; from \eqref{mean2} one gets (now, $x$ is fixed):

\begin{equation} \label{mean3}
E [ \tau  | x_R =u ] =  \frac 1 r \left(1- e^{-x \left (\mu +  \sqrt {\mu ^2 + 2r} \ \right )} \right ) e^{ u \left (\mu +  \sqrt {\mu ^2 + 2r} \ \right ) }.
\end{equation}
Thus, if $h$ is the density of $x_R,$ the mean first-passage time becomes:
\begin{equation}
E[\tau ] = \frac 1 r \left(1- e^{-x \left (\mu +  \sqrt {\mu ^2 + 2r} \ \right )} \right ) \int _0 ^{+ \infty} e^{ u \left (\mu +  \sqrt {\mu ^2 + 2r} \ \right ) } h(u) du .
\end{equation}
Now, we consider the following IMFPT problem:
\begin{equation} \label{IMFPTcaseII}
{given \ m >0, \ find \ the \ density \ h \ of \ x_R \ (if \ it \ exists),\ so \ that \ E[\tau] =m.}
\end{equation}
The {density $h$} is called a solution to the inverse problem \eqref{IMFPTcaseII} (also now, the uniqueness of the solution is not guaranteed).\par
We show below some examples of solutions to the IMFPT  problem \eqref{IMFPTcaseII}. For the sake of simplicity, we limit ourselves to
(undrifted) BM with resetting ($\mu =0);$
to verify them it  suffices to substitute the various quantities into \eqref{mean3}, with $\mu =0.$
\bigskip

\noindent {\bf Example 4.3}
Let $\mathcal X(t)$ be BM with resetting.
For $a, \ \theta >0 ,$ with $\theta > \sqrt {2r},$
if
$$m= \frac {1-e^{-x \sqrt {2r}} } {r } \frac {\theta ^a} {(\theta - \sqrt {2r} ) ^ a },$$
then a solution $h$ to the inverse problem \eqref{IMFPTcaseII} is the Gamma density with parameters $a$ and $\theta,$ that is
$ h(u) = \frac { \theta ^a} { \Gamma (a)} u^{a-1} e^{ - \theta u} \mathbb{I}_{(0, + \infty )}(u).$ \par\noindent
In particular, if $a=1, \ h$ is the exponential density with parameter $\theta.$
\bigskip

\noindent {\bf Example 4.4}
Let $\mathcal X(t)$ be BM with resetting.
For $a, \ \theta >0 ,$ with $\theta > \sqrt {2r},$  if
$$m= \frac {\theta } {r(\theta - \sqrt {2r}) } \ \frac {(1- e^{-x \sqrt {2r}}  ) ( 1- e^{-x ( \theta - \sqrt {2r})} )} { 1- e^{- \theta x}} ,$$
then a solution $h$ to the inverse problem \eqref{IMFPTcaseII} is the truncated exponential density in  the interval $(0,x),$ with parameter $\theta,$ that is $h(u) = \frac {\theta e^{ - \theta u} } {1- e^{- \theta x} } , \ u \in (0,x).$

\bigskip

\noindent {\bf Example 4.5}
Let $\mathcal X(t)$ be BM with resetting.
If $m= \frac 2 {r x \sqrt {2r} }  (\cosh(x \sqrt {2r}) -1),$ then a solution $h$ to the inverse problem \eqref{IMFPTcaseII} is the uniform density in the interval $(0,x).$

\subsection{The Inverse Mean-FET (IMFET) problem}
Let $\mathcal X(t)$ be drifted BM with resetting (with drift $\mu),$ starting from $x \in (0,b),$ where  $x$ and $r$ are fixed, whereas $x_R \in (0,b)$ is random and independent of
$\mathcal X(t). $
Let be
$\tau _{0,b} = \inf \{ t>0: \mathcal X(t) \notin (0,b) \}$ the first-exit time (FET) of $\mathcal X(t)$ from $(0,b);$  then,
$ E[\tau _{0,b}(u)]:= E[\tau _{0,b} |x_R =u] = C_1  (e^{d_1} x -1)+  C_2  (e^{d_2} x -1),$
where  $C_i = C_i (u) $ are given by \eqref{costantsC} and $d_ 1 = -\mu - \sqrt {\mu ^2 + 2r}, \ d_ 2 = -\mu + \sqrt {\mu ^2 + 2r}$ (note that now, unlike the case of Subsection 3.3,  $C_i$
have to be considered as functions of $x_R=u,$ being $r$ and $x$ fixed). Hence, if $h(u), \ u \in (0,b),$ is the density of  $x_R,$
one gets
\begin{equation}
E[ \tau _{0,b} ] = (e^{d_1 x} -1) \int _0 ^b C_1 (u) h(u) du + (e^{d_2 x} -1) \int _0 ^b C_2 (u) h(u) du ,
\end{equation}
and, by using \eqref{costantsC}, we obtain:
\begin{equation} \label{m3}
 E[ \tau _{0,b} ] = \frac 1 r \left [  (e^{d_1x } -1) (1- e^{d_2b}) - (e^{d_2x } -1) (1- e^{d_1b}) \right ] \times
\end{equation}
$$
\times  \int _0 ^b \frac {e^{-d_2 u} h(u) } { 1- e^{d_1b} -e^{-u(d_1-d_2)} (1-e^{d_2b})} du.
$$
Now, one can consider the following kind of inverse problem for the mean of the FET (IMFET) problem):
\begin{equation} \label{IMFETcaseII}
{given \ m >0, \ find \ the \ density \ h \ of \ x_R \ (if \ it \ exists),\ so \ that \ E[\tau _{0,b}] =m.}
\end{equation}
The {density $h$}  is called a solution to the inverse problem \eqref{IMFETcaseII}; in fact, also in this case the uniqueness of solution is not guaranteed.
\bigskip

\noindent {\bf Example 4.6}
Let $\mathcal X(t)$ be BM with resetting
and let
$$ m= \frac 1 {r \sqrt {2r}  (1-e^{-b \sqrt {2r}})} \left [  (e^{ -x \sqrt {2r} } -1) (1- e^{b \sqrt {2r}}) - (e^{x \sqrt {2r} } -1) (1- e^{- b \sqrt {2r}}) \right ] \times $$
$$ \times  \left [ 1- e^ {-b \sqrt {2r}} + \frac { \sqrt {\widetilde D} } 2  \left ( \ln  \frac {1- \sqrt {\widetilde D }} { 1+ \sqrt {\widetilde D } }   - \ln  \frac {e^ {-b \sqrt {2r}}- \sqrt {\widetilde D }} { e^{-b \sqrt {2r}}+ \sqrt {\widetilde D}  } \right ) \right ]$$
where $\widetilde D= \frac {1-e^{b \sqrt {2r}} } {1-e^{-b \sqrt {2r}} }  ;$
then, the uniform density $h$ on the interval $(0,b)$ is a solution to the IMFET problem \eqref{IMFETcaseII}. \par\noindent
This soon follows by substituting the various quantities in \eqref{m3}.

\section{Conclusions and Final Remarks}
{We studied several types of inverse problems for the first-passage place and the first-passage time of a one-dimensional diffusion process $\mathcal X(t)$ with resetting,
obtained from an underlying
temporally homogeneous diffusion process $X(t),$
driven by the SDE
$
dX(t)=  \mu(X(t)) dt + \sigma (X(t)) d {W_t} ,
$
and starting from $\mathcal X(0)= X(0)=\eta ,$
where ${W_t}$ is standard Brownian motion, and the drift $\mu (\cdot)$ and diffusion coefficient $\sigma (\cdot)$ are regular enough functions, such that there exists a unique strong solution of the SDE.
}
\par
{Actually,
we considered two cases; in the first one the initial position
$\eta = \mathcal X(0)$
 was random  and independent of $\mathcal X(t),$  whereas the resetting rate $r$ and the reset position $x_R$ were fixed; in the second case $\eta = \mathcal X(0)$ and $r$ were fixed, whereas  $x_R$ was random and independent of $\mathcal X(t) .$ }\par
{In general, the solution to the inverse first-passage place (IFPP) problem is not unique; instead, the solution to the inverse first-passage time (IFPT) problem is unique, provided that its Laplace transform is analytic.
}
\par
{For all the  inverse problems, we reported several explicit examples of solutions, mostly concerning Wiener process with resetting.
We did not deal with the analogous inverse problem  for $\mathcal X(t),$ which consists in finding the density of $\eta \in (0,b)$ (or of the reset position $x_R \in (0,b)),$ in such a way that the first-exit time of  $\mathcal X(t)$ from $(0,b)$ has an assigned distribution; this, as well as the corresponding direct problem will be the subject of a future investigation.
}
\par
{We remark that the inverse problems considered concern randomization in the starting point $\eta $ of $\mathcal X(t),$ or in the reset position $x_R;$
more generally, one could introduce randomization both in the starting point $\eta$ and in the reset position $x_R,$ and then
study the corresponding inverse problems, where now a solution is the joint density of $(\eta, \ x_R),$ if it exist;
this too will possibly be the subject of a future article.
}
\par
{Our study was motivated by the fact that, as in the case without resetting, direct and inverse problems for first-passage time and first-passage place of a diffusion with resetting
$\mathcal X(t)$ are very worthy of attention, because they have notable applications in several applied fields, e.g.
in biological modeling concerning neuronal activity,  queuing theory, and  mathematical finance.
}

\section*{Acknowledgments}
\noindent
{The author expresses particular thanks to the anonymous Reviewers for their
valuable comments, leading to an improved presentation.
The author belongs to GNAMPA, the Italian National Research Group of INdAM; he also acknowledges the MIUR Excellence Department Project MatMod@TOV awarded to the Department of Mathematics, University of Rome Tor Vergata, CUP E83C23000330006. }
\bigskip

\bibliographystyle{amsplain}

\end{document}